\documentclass{agtart_a}
\pdfoutput=1
\usepackage[all]{xy}

\title[A lower bound for coherences on $\BP$]{A lower bound for coherences on the\\Brown--Peterson spectrum}  
\author{Birgit Richter}
\givenname{Birgit}
\surname{Richter}
\address{Fachbereich Mathematik der Universit\"at Hamburg\\
Bundesstra{\ss}e 55\\\newline 20146 Hamburg\\Germany}
\email{richter@math.uni-hamburg.de}
\urladdr{http://www.math.uni-hamburg.de/home/richter/}


\received{25 May 2005}
\revised{17 November 2005}
\accepted{14 February 2006}
\proposed{}
\seconded{}
\publishedonline{26 February 2006}
\published{26 February 2006}

\volumenumber{6}
\issuenumber{}
\publicationyear{2006}
\papernumber{9}
\startpage{287}
\endpage{308}

\doi{}
\MR{}
\Zbl{}

\arxivreference{math.AT/0504322}
\arxivpassword{}

\let\xysavmatrix\xymatrix
\def\xymatrix{\disablesubscriptcorrection\xysavmatrix}
\AtBeginDocument{\let\bar\wbar}
\newcommand{\we}{\smash{\rlap{\kern 6pt 
\raise 4pt\hbox{\footnotesize $\sim$}}}\longrightarrow}

\makeatletter
\def\cnewtheorem#1[#2]#3{\newtheorem{#1}{#3}[section]
\expandafter\let\csname c@#1\endcsname\c@thm}

\newtheorem{thm}{Theorem}[section]
\cnewtheorem{lem}[thm]{Lemma}
\cnewtheorem{prop}[thm]{Proposition}
\cnewtheorem{cor}[thm]{Corollary}

\theoremstyle{definition}
\cnewtheorem{rem}[thm]{Remark}
\cnewtheorem{defn}[thm]{Definition}
\cnewtheorem{example}[thm]{Example}
\cnewtheorem{conj}[thm]{Conjecture}

\makeatother  %  move after \newtheorem block

\def\HG{\mathrm{H}\Gamma}

\def\thhbp{\mathit{THH}(BP)}
\def\BSF{\mathit{BSF}}
\def\SF{\mathit{SF}}
\def\BP{\mathit{BP}}
\def\BSO{\mathit{BSO}}
\def\Z{\mathbb{Z}}
\def\F{\mathbb{F}}
\def\lra{\longrightarrow}
\def\leq{\leqslant}
\def\geq{\geqslant}
\def\lie{\mathrm{Lie}}
\def\ra{\rightarrow}
\DeclareMathOperator{\colim}{colim}

\def\bpi{\BP\langle i\rangle}
\def\bp1{\BP\langle 1\rangle}
\def\ie{ie}
\begin{document}

\begin{asciiabstract}
We provide a lower bound for the coherence of the homotopy
commutativity of the Brown-Peterson spectrum, BP, at a given prime p
and prove that it is at least (2p^2 + 2p - 2)-homotopy commutative. We
give a proof based on Dyer-Lashof operations that BP cannot be a Thom
spectrum associated to n-fold loop maps to BSF for n=4 at 2 and n=2p+4
at odd primes.  Other examples where we obtain estimates for coherence
are the Johnson-Wilson spectra, localized away from the maximal ideal
and unlocalized. We close with a negative result on Morava-K-theory.
\end{asciiabstract}

\begin{htmlabstract}
We provide a lower bound for the coherence of the homotopy commutativity
of the Brown&ndash;Peterson spectrum <i>BP</i> at a given prime <i>p</i>
and prove that it is at least <i>(2p&sup2;+2p-2)</i>&ndash;homotopy
commutative. We give a proof based on Dyer&ndash;Lashof operations that
<i>BP</i> cannot be a Thom spectrum associated to <i>n</i>&ndash;fold
loop maps to <i>BSF</i> for <i>n=4</i> at 2 and <i>n=2p+4</i> at odd
primes.  Other examples where we obtain estimates for coherence are the
Johnson&ndash;Wilson spectra, localized away from the maximal ideal and
unlocalized. We close with a negative result on Morava-K&ndash;theory.
\end{htmlabstract}

\begin{abstract}
We provide a lower bound for the coherence of the homotopy
commutativity of the Brown--Peterson spectrum, $BP$, at a given prime
$p$ and prove that it is at least $(2p^2 + 2p - 2)$--homotopy
commutative. We give a proof based on Dyer--Lashof operations that $BP$
cannot be a Thom spectrum associated to $n$--fold loop maps to $BSF$
for $n=4$ at $2$ and $n=2p+4$ at odd primes.  Other examples where we
obtain estimates for coherence are the Johnson--Wilson spectra,
localized away from the maximal ideal and unlocalized. We close with a
negative result on Morava-$K$--theory.
\end{abstract}

\maketitle
\section{Introduction}

Recently obstruction theory for imposing $E_\infty$--structures on
homotopy commutative and associative ring spectra has been
successfully applied in some cases of very well-behaved higher
chromatic spectra like $E_n$ (Goerss--Hopkins \cite{GH},
Richter--Robinson \cite{RRo}) and $\widehat{E(n)}$ (Baker--Richter
\cite{BR}) and up to chromatic type one \cite{BR}.

The two approaches for such an obstruction theory that are available
at the moment were developed by Paul Goerss and Mike Hopkins \cite{GH}
on the one hand and Alan Robinson \cite{Ro:obstr} on the other
hand. They use Andr\'e--Quillen type cohomology theories as the habitat
for possible obstructions which turn out to be equivalent (see
Basterra--Richter \cite[Theorem 2.6]{BaR}).  The crucial point in all
the examples mentioned above is that one can rely on an \'etaleness
property of the algebra of cooperations in order to make all
obstruction groups vanish.

It is an old, still open question raised by Peter May, whether the
Brown--Peterson spectrum, $\BP$, is an $E_\infty$ ring spectrum. There
have been serious attempts to prove that $\BP$ has an $E_\infty$
model. In fact, topological Andr\'e--Quillen homology defined by Maria
Basterra in \cite{Ba} was originally introduced to solve this problem
(compare Kriz \cite{Kriz}).

Unlike
complex cobordism, $\mathit{MU}$, or other Thom spectra, $\BP$ is not born with
an $E_\infty$--structure. 
Therefore trying obstruction theory methods seems to be a canonical attempt. 
But the  algebra of cooperations $\BP_*\BP$, is a polynomial algebra on 
countably many generators, and this will produce a lot of non-trivial 
Andr\'e--Quillen type cohomology groups.  

However, the possible obstruction groups for instance in the setting
of Robinson's Gamma cohomology only occur in a small range of degrees. The
approach of this note is to exploit the sparseness of these degrees
together with the sparseness of the coefficients 
$$\BP_* = \Z_{(p)}[v_1,v_2,\ldots], \text{ with } |v_i| = 2p^i -2$$
and the algebra $\BP_*\BP$ to obtain an estimate for the coherence of the 
homotopy commutativity of $\BP$. We admit that this lower bound might be much
too pessimistic, but it seems that there is some interest in even
partial results. 

Maria Basterra and Mike Mandell announced that $\BP$ 
has at least an $E_4$--structure, ie, possesses an action by
the  little $4$--cubes operad. By results of Fiedorowicz--Vogt \cite{FV}
or Basterra--Mandell this implies that the topological
Hochschild homology spectrum of $\BP$, $\thhbp$,  is at
least an $E_3$--spectrum and the natural map $\BP \ra \thhbp$
is a map of $E_3$--spectra. 

Our approximation result depends on the prime involved. 
\begin{thm}\label{thm:bp}
The Brown--Peterson spectrum $\BP$ at a prime $p$ has
at least a $(2p^2 + 2p -2)$--stage structure.
\end{thm}

Here, an $n$--stage structure is a certain filtration step towards an 
$E_\infty$--structure. Alan Robinson introduced this filtration in 
\cite{Ro:obstr} and we will describe it in 
\fullref{filtration}. We will prove \fullref{thm:bp} in
\fullref{subsec:theproof}.

The obstruction theory which we will use is the one developed by Alan
Robinson \cite{Ro:obstr,Ro:units}. We will explain it in detail in
sections \ref{sec:gammahom}--\ref{sec:geom}. 

At the moment we do not have a direct comparison between  $n$--stage
structures and structures over some version of the little cubes
operad. 
% Our  conjecture is that an $(n+2)$--stage structure which is
% compatible  with a given $A_\infty$--structure gives rise to an 
% $E_n$--structure. 

However, $n$--stage structures are of independent interest because they
give rise to  Dyer--Lashof operations in a certain range of degrees 
(see 
\fullref{sec:dl}). In the case of $\BP$ one could get homology operations by
comparing the $\F_p$--homology of $\BP$ with the one of $\mathit{MU}$
and  $H\F_p$
(see Bruner--May--McClure--Steinberger \cite[page 63]{Hinfty}). 
%%%changed
% But as we can prove that some low dimensional 
% operations have their geometric origin in $\BP$ itself,
We  use Dyer--Lashof operations 
to show in \fullref{thm:nothom} that $\BP$ cannot be the Thom spectrum 
%%%changed
associated to an $n$--fold loop map to $\BSF$. Here $n$ is four for the
even prime and $2p+4$ for any odd prime $p$. We stress that a stronger
result is mentioned in Lewis' thesis \cite{Lewis}.

In other examples the existence of Dyer--Lashof operations might
help for instance with calculations of topological
Hochschild homology. 

The main object of study of this note is the Brown--Peterson
spectrum. The reader who is primarly interested in this example might
skip sections \ref{sec:gammahom}--\ref{sec:geom} and might procede
directly to \fullref{sec:bp}. The material presented in sections
\ref{sec:gammahom}--\ref{sec:geom} is needed for cases like the 
Johnson--Wilson spectra. We discuss this and other examples in
\fullref{sec:examples}. \par \bigskip \noindent
{\bf Acknowledgements} \, I thank Sverre Lun{\o}e-Nilsen and John
Rognes for their interest. They pointed out to me that my partial
coherence results should lead to Dyer--Lashof operations. Thanks to
Andy Baker who read earlier versions of this and made important
comments. I am grateful to John Rognes and the Department of
Mathematics in Oslo for their hospitality. Parts of this work were written
when I was in Bonn and I would like to thank Carl-Friedrich
B\"odigheimer for his constant support during the last years. The
author was  partially supported by
the \emph{Strategisk Universitetsprogram i Ren Matematikk} (SUPREMA) 
of the Norwegian Research Council. 
\section{Some background on Gamma (co)homology} \label{sec:gammahom}
Let $k$ be a (graded) commutative ring with unit, let $A$ be a (graded)
commutative $k$--algebra and let $M$ be a (graded) $A$--module. In the
following tensor products will all be taken with respect to $k$. Gamma
homology of $A$ over $k$ with coefficients in $M$ is defined in 
Robinson \cite[2.5]{Ro:obstr} as the homology of the total complex of a
bicomplex $\Xi_{*,*}$ which we will now describe.

Let $\lie(n)$ be the $n^{\rm th}$ term of the operad which codifies
Lie algebras over $k$, ie, $\lie(n)$ is the free $k$--module
generated by all Lie monomials in variables $x_1,\ldots,x_n$ such that
each variable appears exactly once. There is a canonical action of the
symmetric group on $n$ letters, $\Sigma_n$, on $\lie(n)$ by permuting
the variables $x_i$. Let $\lie(n)^*$ be the $k$--linear dual of
$\lie(n)$. Then the bicomplex for Gamma homology in bidegree $(r,s)$
is defined as
$$\Xi_{r,s}(A|k;M) = \lie(s+1)^* \otimes k[\Sigma_{s+1}]^{\otimes r} 
\otimes A^{\otimes (s+1)} \otimes M. $$
Here, all undecorated tensor products are taken over the ground ring $k$. 
The horizontal differential is the differential of the bar construction;  
% ie, the elements of the symmetric group are multiplied
% together or an action of $\Sigma_{s+1}$ on the dual of the Lie
% monomials is induced or the elements in the $(s+1)^{\rm st}$ tensor power of 
% the algebra $A$ are permuted. 
the vertical differential is complicated  and
we refer the curious reader to \cite[Secion 2]{Ro:obstr} for details. For
our purpose it is enough to know that it induces multiplication of the
algebra entries, induces an action of $A$ on $M$ and does something
with the permutations in $\Sigma_{s+1}$ to reduce them to elements in
$\Sigma_{s}$.

We emphasize that all these operations preserve the internal degree if
$A$ and $M$ are graded. Therefore it makes sense to define the $i^{\rm th}$
homogeneous part of $\Xi_{*,*}$ and the associated total
complex. Gamma homology for graded commutative algebras therefore
possesses a natural bigrading, where $\HG_{q,i}(A|k;M)$ is the $q^{\rm th}$
homology of the $i^{\rm th}$ homogeneous part of the total complex of
$\Xi_{*,*}(A|k;M)$. Gamma cohomology, which is defined via the
homomorphism complex out of the total complex of $\Xi_{*,*}(A|k;A)$
into $M$, inherits the internal grading. Following \cite{Ro:obstr} we
denote by $\HG^{q,i}(A|k;M)$ the $q^{\rm th}$ cohomology of the homomorphism
complex 
$$ \mathrm{Hom}_{A}^i(\mathrm{Tot}(\Xi_{*,*}(A|k;A)), M)$$
whose morphisms lower degree by $i$. 

\section{Robinson's obstruction theory} \label{sec:obs}
Consider the product of the topological version of the Barratt--Eccles
operad 
$(E\Sigma_n)_n$ and Boardman's tree operad $(T_n)_n$. Here
the space of $n$--trees, $T_n$,  consists of abstract trees on $n+1$
leaves. These leaves are labelled with the numbers $0,\ldots,n$ where
each label appears exactly once. Internal edges get an assigned length
$0 < \lambda \leq 1$. This
tree space is set to consist of a point for $n \leq 2$, the only
$2$--tree is the tree
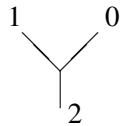
\vspace{0.3cm}

\begin{center} 
\begin{picture}(4,4)
\setlength{\unitlength}{1cm}
\put(0,0){\line(1,1){0.5}}
\put(0.6,0.6){0}
\put(0,0){\line(-1,1){0.5}}
\put(-0.7,0.6){1}
\put(0,0){\line(0,-1){0.5}}
\put(0.1,-0.7){2}
\end{picture}
\end{center}
%\mbox{}\\[0.3cm]
\vspace{0.3cm}

There are three different types of $3$-trees, namely 
%\mbox{}\\[0.3cm]
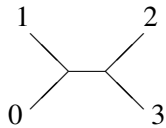
\vspace{0.3cm}

\begin{center} 
\begin{picture}(4,4)
\setlength{\unitlength}{1cm}
\put(0,0){\line(-1,-1){0.5}}
\put(-0.8,-0.7){0}
\put(0,0){\line(-1,1){0.5}}
\put(-0.7,0.6){1}
\put(0,0){\line(1,0){0.5}}
\put(0.5,0){\line(1,1){0.5}}
\put(1,0.6){2}
\put(0.5,0){\line(1,-1){0.5}}
\put(1.1,-0.7){3}
\end{picture}
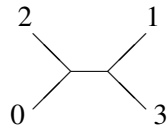
\hspace{4cm}
\begin{picture}(4,4)(4,0)
\setlength{\unitlength}{1cm}
\put(0,0){\line(-1,-1){0.5}}
\put(-0.8,-0.7){0}
\put(0,0){\line(-1,1){0.5}}
\put(-0.7,0.6){2}
\put(0,0){\line(1,0){0.5}}
\put(0.5,0){\line(1,1){0.5}}
\put(1,0.6){1}
\put(0.5,0){\line(1,-1){0.5}}
\put(1.1,-0.7){3}
\end{picture}
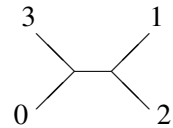
\hspace{4cm}
\begin{picture}(4,4)
\setlength{\unitlength}{1cm}
\put(0,0){\line(-1,-1){0.5}}
\put(-0.8,-0.7){0}
\put(0,0){\line(-1,1){0.5}}
\put(-0.7,0.6){3}
\put(0,0){\line(1,0){0.5}}
\put(0.5,0){\line(1,1){0.5}}
\put(1,0.6){1}
\put(0.5,0){\line(1,-1){0.5}}
\put(1.1,-0.7){2}
\end{picture}
\end{center}
%\mbox{}\\[0.3cm]
\vspace{0.3cm}

with a corolla-shaped tree if the length of the only internal edge is
shrunk to zero.
%\mbox{}\\[0.3cm]
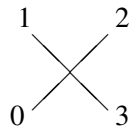
\vspace{0.3cm}

\begin{center} 
\begin{picture}(4,4)
\setlength{\unitlength}{1cm}
\put(0,0){\line(-1,-1){0.5}}
\put(-0.8,-0.7){0}
\put(0,0){\line(-1,1){0.5}}
\put(-0.7,0.6){1}
\put(0,0){\line(1,1){0.5}}
\put(0.6,0.6){2}
\put(0,0){\line(1,-1){0.5}}
\put(0.6,-0.7){3}
\end{picture}
\end{center}
%\mbox{}\\[0.3cm]
\vspace{0.3cm}

For arbitrary $n$, the space $T_n$ is contractible with the corolla on
$n+1$ leaves as basepoint. Composition in the tree operad is given by
grafting trees. The newly built internal edge in the composed tree  is
defined to be of length one. The  $n^{\rm th}$ part of the  operad which
Robinson uses  is the product operad 
$$ \mathcal{B}_n = E\Sigma_n \times T_n.$$
As $E\Sigma_n$ is $\Sigma_n$--free and contractible and as $T_n$ is
contractible, the product operad $\mathcal{B}$ is an $E_\infty$--operad. 
Robinson defines a filtration of this $E_\infty$--operad as follows:
set $\mathcal{B}_n^{(i)} := (E\Sigma_n)^{(i)} \times T_n$
where $(E\Sigma_n)^{(i)}$ is the $i^{\rm th}$ skeleton of the standard model
for $E\Sigma_n$. Then define 
\begin{equation}
\nabla^n \mathcal{B}_m := \mathcal{B}_m^{(n-m)}.
\end{equation}
\begin{defn}\cite[5.3]{Ro:obstr}\qua
An $n$--stage structure for an $E_\infty$--structure on a spectrum $E$
is a sequence of maps
$$ \mu_m\co \nabla^n\mathcal{B}_m \ltimes_{\Sigma_m} E^{\wedge m} \lra E$$
which on their restricted domain of definition satisfy the
requirements for an operad action on $E$.
\end{defn}
Let us make explicit what that amounts to in small filtration
degrees. \label{filtration}
A $2$--stage structure on a spectrum $E$ consists of action maps
starting from $\nabla^2\mathcal{B}_m$ which is
$$\mathcal{B}_m^{(2-m)} = (E\Sigma_m)^{(2-m)} \times T_m.$$ 
The $(2-m)$--skeleton of $E\Sigma_m$ is trivial for $m > 2$, therefore
the only  requirement is that we have a map 
$$ ((E\Sigma_1)^{(1)} \times {T}_1) \ltimes_{\Sigma_1} E 
\cong E  \stackrel{\varphi}{\lra} E $$
and that $E$ possesses a map
$$ ((E\Sigma_2)^{(0)} \times {T}_2) \ltimes_{\Sigma_2}
E^{\wedge 2} \cong (\Sigma_2) \ltimes_{\Sigma_2} E^{\wedge 2} \lra E. $$
This is nothing but a multiplication $\mu$ on $E$ together with its
twisted version $\mu \circ \tau$ if $\tau$ denotes the generator of
$\Sigma_2$. The axioms of an operad action force the map $\varphi$ to
be the identity of $E$. Iterates of $\mu$ and $\mu \circ \tau$ act on
higher smash powers of $E$, but they do not have to satisfy any
relations. 

A $3$--stage structure on $E$ comes with three kinds of maps,  because 
non-trivial values for $m$ are $1,2,3$. The first value does not give
anything new. The second step requires action maps 
$$ ((E\Sigma_2)^{(1)} \times {T}_2) \ltimes_{\Sigma_2}
E^{\wedge 2} \lra E. $$
The $1$--skeleton of $E\Sigma_2$ is the $1$--circle, giving the homotopy
between $\mu$ and $\mu \circ \tau$.

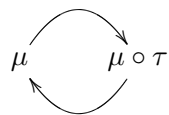
\smallskip

$$
\xymatrix{
{\mu}\ar@/^4ex/[r] & {\mu \circ \tau}\ar@/^4ex/[l]
}
$$

\medskip

\smallskip
\noindent
In addition to that, the value $m=3$ brings in the homotopies for
associativity via the trees we described above, because we obtain maps 
$$ ((E\Sigma_3)^{(0)} \times {T}_3) \ltimes_{\Sigma_3}
E^{\wedge 3} \cong (\Sigma_3 \times T_3) \ltimes_{\Sigma_3}
E^{\wedge 3} \lra E. $$
\begin{thm}{\rm \cite[Theorem 5.5]{Ro:obstr}}\qua
Assume that $E$ is  a homotopy commutative and associative ring
spectrum which satisfies 
\begin{equation} \label{uckuenneth}
E^*(E^{\wedge m}) \cong \mathrm{Hom}_{E_*}(E_*E^{\otimes m}, E_*)
\quad \mathrm{for \phantom{q}  all} \quad m \geq 1.
\end{equation}
If $E$ has an $(n-1)$--stage structure which can be extended to an
$n$--stage structure then possible obstructions to
extending this further to an $(n+1)$--stage structure live in 
$$ \HG^{n,2-n}(E_*E | E_*; E_*). $$
If in addition $\HG^{n,1-n}(E_*E | E_*; E_*)$ vanishes, then this
extension is unique. 
\end{thm}
We start with a $3$--stage structure. If we want to establish an
$n$--stage structure on $E$, then we have to show that Gamma
cohomology vanishes in bidegrees $(\ell,2-\ell)$ for all 
$n-1 \geq \ell \geq 3$. (This is incorrectly stated in \cite[5.6]{Ro:obstr}
but tacitly corrected in \cite[5.8]{Ro:units}.)

Let us note that $\BP$ satisfies the necessary properties to apply
Robinson's obstruction theory: $\BP$ is a homotopy commutative
$\mathit{MU}$--ring spectrum at all primes, Strickland
\cite[2.8,2.9]{St}, and it satisfies the requirement from
\ref{uckuenneth}.

\section{Where the obstructions come from} \label{sec:geom}
In this section we will explain the geometric origin 
of Robinson's obstruction groups \cite{Ro:obstr}. We do not
present anything new here, but focus on some details that we will need
later in \fullref{sec:examples}. 

The obstruction groups in the theorem above arise from a geometric
object: one can extend an $n$--stage structure to an $(n+1)$--stage
structure if certain cohomology classes vanish. An $n$--stage
always gives rise to maps out of parts of  $\nabla^{n+1}
\mathcal{B}_m$ as well, because we require that the action maps
$\mu_m$ satisfy the axioms of operad actions where they are
defined. Therefore these actions are closed under
composition. Robinson calls the part in $\nabla^{n+1} \mathcal{B}_m$
that is generated by compositions the \emph{boundary} of $\nabla^{n+1}
\mathcal{B}_m$ and denotes it by $\partial \nabla^{n+1}
\mathcal{B}_m$. Let $Q_{n+1}^m$ be the cofibre of the
inclusion of $\partial \nabla^{n+1}
\mathcal{B}_m \cup \nabla^n \mathcal{B}_m$ into $\nabla^{n+1}
\mathcal{B}_m$ and  consider the following diagram: 
$$
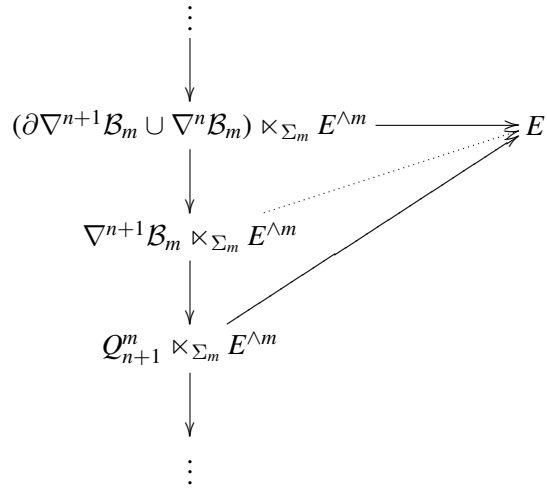
\xymatrix{
{\vdots}\ar[d]  & &   \\
{(\partial \nabla^{n+1}
\mathcal{B}_m \cup \nabla^n \mathcal{B}_m) \ltimes_{\Sigma_m} E^{\wedge 
m}} \ar[d] \ar[rr] & & E\\
{\nabla^{n+1} \mathcal{B}_m \ltimes_{\Sigma_m} E^{\wedge m}} \ar[d]
\ar@{.>}[rru] && \\
{Q_{n+1}^m \ltimes_{\Sigma_m} E^{\wedge m}}\ar[d]
\ar[rruu] && \\
{\vdots} & & 
}
$$
We obtain a long exact sequence in $E$--cohomology 
\begin{multline*}
 \ldots \ra E^0(Q_{n+1}^m \ltimes_{\Sigma_m} E^{\wedge m}) \ra 
E^0(\nabla^{n+1} \mathcal{B}_m \ltimes_{\Sigma_m} E^{\wedge m}) \ra \\
\ra E^0((\partial \nabla^{n+1}
\mathcal{B}_m \cup \nabla^n \mathcal{B}_m) \ltimes_{\Sigma_m}
E^{\wedge m}) \ra E^1(Q_{n+1}^m \ltimes_{\Sigma_m} E^{\wedge m}) \ra \ldots
%\phantom{\nabla^{n+1} \mathcal{B}_m \ltimes_{\Sigma_m} E^{\wedge m}) \ra 
% E^0((\partial \nabla^{n+1}
% \mathcal{B}_m \cup \nabla^n \mathcal{B}_m) \ltimes_{\Sigma_m}
% E^{\wedge m}) \ra}
\end{multline*}
If the first $E$--cohomology group of 
$Q_{n+1}^m \ltimes_{\Sigma_m} E^{\wedge m}$ vanishes, then one can 
extend the action map to the
$(n+1)^{\rm st}$ filtration step. If the term $E^0(Q_{n+1}^m
\ltimes_{\Sigma_m} E^{\wedge m})$  is trivial as well, then the
extension of the $n$--stage structure to an $(n+1)$--stage structure is
unique. Note, that we have to consider all $2 \leq m \leq n+1$, because
these are the non-trivial domains for an $(n+1)$--stage action. 

If one has a universal coefficient theorem  at
hand, one can then continue to identify this cohomology group as 
$\mathrm{Hom}_{E_*}(E_*(Q_{n+1}^m \ltimes_{\Sigma_m} E^{\wedge m}),
E_*)$. Note, that up to this stage homomorphisms are taken with
respect to $E_*$ and not $E_*E$. Alan Robinson then identifies 
$E_*(Q_{n+1}^m \ltimes_{\Sigma_m} E^{\wedge m})$ as some tractable
algebraic object and induces
homomorphism up to $E_*E$. Furthermore he proves in \cite[Proposition
5.4]{Ro:obstr} that obstruction 
are always cocycles for a suitable coboundary map, which then gives the
final identification with Gamma cohomology groups. 

However, one could stop at earlier stages before passing to the level
of Gamma cohomology groups. One can try to extract direct information
either out of 
\begin{equation} \label{cohomology_level}
  E^1(Q_{n+1}^m \ltimes_{\Sigma_m} E^{\wedge m})
\end{equation}
or, in the presence of a universal coefficient theorem for $E$, one
can investigate in which degrees the groups 
\begin{equation} \label{homology_level}
  \mathrm{Hom}_{E_*}^{1}(E_*(Q_{n+1}^m \ltimes_{\Sigma_m} E^{\wedge m}),
E_*) 
\end{equation}
are non-trivial. 

In the examples that we will consider, we already know that the
spectra under consideration possess $A_\infty$--structures. Therefore
we have universal coefficient spectral sequences available converging
to the cohomological term 
\ref{cohomology_level}. 

Let us examine the $E$--homology of $Q_{n+1}^m \ltimes_{\Sigma_m}
E^{\wedge m}$ further in order to get control over internal degree
shifts. 

We use Robinson's geometric identification of $Q_{n+1}^m$ as
$E\Sigma_m^{(n-m+1)}/E\Sigma_m^{(n-m)} \wedge {T}_m/\partial
T_m$. Here the second quotient is the quotient of the cubical tree
complex modulo the subcomplex of all fully grown trees -- these are
decomposable trees, \ie, trees with at least one internal edge having
length one. The first quotient is equivalent to
$\bigvee_{\Sigma_m^{n-m+2}} \mathbb{S}^{n-m+1}$ whereas the
identification in \cite{Ro:obstr} and Whitehouse \cite{Wh} shows that
the second quotient has the homotopy type of $\bigvee_{(m-1)!}
\mathbb{S}^{m-2}$ which on the level of homology gives rise to the
dual of the Lie representation. One copy of $\Sigma_m$ is swallowed by
the $\ltimes_{\Sigma_m}$. The degree shift caused by spheres adds up
to a total of $n-1$, so in $E$--homology we get the term
$$ \Sigma^{n-1} \lie(m)^* \otimes E_*[\Sigma_m]^{\otimes (n-m+1)} \otimes
E_*E^{\otimes m}.$$
Therefore the $E_*$--homomorphisms on these are given as
\begin{multline*}
\mathrm{Hom}_{E_*}^{1}(E_*(Q_{n+1}^m \ltimes_{\Sigma_m} E^{\wedge m}),
E_*) \\
\cong \mathrm{Hom}_{E_*}^{2-n}(\lie(m)^* \otimes E_*[\Sigma_m]^{\otimes
(n-m+1)} \otimes E_*(E^{\wedge m}), E_*).  
\end{multline*}

\section{The case of $\BP$} \label{sec:bp}
Potential obstruction classes for an extension to an $(n+1)$--stage
structure  on $\BP$ live in bidegree $(n,2-n)$. Here $n$ is the cohomological 
degree and $2-n$ is the internal degree. 
The algebra of cooperations $\BP_*\BP$ over $\BP_*$ is a polynomial ring
over $\BP_*$, 
$$ \BP_*\BP = \BP_*[t_1,t_2,\ldots] \text{ with } |t_i| = 2p^i-2.$$
We are therefore in the situation where we can apply the
universal coefficient theorem to see that
the obstruction groups arise from  
\begin{multline*}
\mathrm{Hom}_{\BP_*}^{1}(\BP_*(Q_{n+1}^m \ltimes_{\Sigma_m} \BP^{\wedge
  m}), \BP_*) \\
\cong \mathrm{Hom}_{\BP_*}^{2-n}(\lie(m)^* \otimes 
\BP_*[\Sigma_m]^{\otimes (n-m+1)} \otimes \BP_*\BP^{\otimes m}, \BP_*).
\end{multline*}
In order to achieve an internal degree of the form
$2-n$, we have to look for a map that raises degree by $n-2$. 
%%%[\Sigma^n X, Y] = [X,Y]_n \ra Hom^n_{E_*E}(E_*X, E_*Y)
Everything in sight is concentrated in degrees of the form 
$$  i = \sum_{i=1}^N \lambda_j (2p^j -2) $$
where the $\lambda_j$ are non-negative integers. 
Any non-trivial map 
$$\lie(m)^* \otimes \BP_*[\Sigma_m]^{\otimes (n-m+1)} \otimes
\BP_*\BP^{\otimes m} \ra \BP_*$$
can only alter the  degree again by a degree of the form $i = \sum
\lambda_j (2p^j -2)$, because $\BP_*$ is concentrated in the same
degrees. The  minimal such degree bigger than zero is $2p-2$. 

Obstruction groups occur as Gamma cohomology of
bidegree $(n,2-n)$ for $n \geq 3$. Therefore the first possibly
non-trivial  obstruction group could live in $\HG^{2p,2-2p}$, which
could contain an obstruction to extending a $2p$--stage structure to a
$(2p+1)$--stage structure. This amounts to saying that $\BP$ has at
least a $2p$--stage structure.

\section{Refined estimates} \label{sec:bpbetter}
The previous section dealt with an argument which was merely a degree
count. We will improve the $2p$--estimate by having a detailed look
at the actual homological level. 
We consider the universal coefficient spectral sequence for Gamma
cohomology of $\BP_*\BP$. Its $E_2$--term  is of the form 
\begin{equation} \label{UCspsq} 
\mathrm{Ext}^{*,*}_{\BP_*\BP}(\HG_{*,*}(\BP_*\BP|\BP_*;\BP_*\BP), \BP_*) 
\Rightarrow \HG^{*,*}(\BP_*\BP|\BP_*; \BP_*). 
\end{equation}

\begin{rem}
In the spectral sequence \ref{UCspsq} the cohomological degree $s$ for
our possible obstruction group $\HG^{s,2-s}$ is now spread over the
groups 
\begin{multline*}
\mathrm{Ext}^{0,s}_{\BP_*\BP}(\HG_{*,*}(\BP_*\BP|\BP_*;\BP_*\BP), 
\BP_*),\\
  \ldots, \mathrm{Ext}^{s,0}_{\BP_*\BP}(\HG_{*,*}(\BP_*\BP|\BP_*;\BP_*\BP),
\BP_*)
\end{multline*}
so we have to prove, that all these groups vanish. As we grade
homologically the
internal degree $j$ in $\mathrm{Ext}^{i,j}$ corresponds to a map that
lowers degree by $j$. For instance  a homomorphism 
$f$ in  $\mathrm{Hom}^j_{\BP_*\BP}(\HG_{*,*}(\BP_*\BP|\BP_*;\BP_*\BP), \BP_*)$ is a
sequence of maps 
$$ f = (f_s); \quad f_s\co \HG_{s,*}(\BP_*\BP|\BP_*;\BP_*\BP) \lra \BP_{s-j}.$$ 
\end{rem}

The reformulations in the rest of this section will not preserve the
internal grading of Gamma homology. The following result is an
immediate application of Robinson--Whitehouse \cite[6.8(2)]{RoWh}.
\begin{lem}\label{lem:splitting}
Gamma homology groups for $\BP_*\BP = \BP_*[t_1,t_2,\ldots]$ split into\break
Gamma homology groups for the single pieces $\BP_*[t_i]$:
\begin{equation} \label{splitting}
\HG_{s,*}(\BP_*\BP|\BP_*;\BP_*\BP) \cong \bigoplus_{i \geq 1}
\HG_{s,*}(\BP_*[t_i]|\BP_*;\BP_*\BP). 
\end{equation}
\end{lem}

We start with simplifying each summand in the splitting
\ref{splitting}.
\begin{lem}
Gamma homology of $\BP_*[t_i]$ can be expressed as follows:
\begin{align*}\HG_{s,*}(\BP_*[t_i]|\BP_*;& \BP_*\BP)\\ &\cong
\BP_* \otimes_{\Z_{(p)}} 
\HG_{s,*}(\Z_{(p)}[t_i]|\Z_{(p)}; \Z_{(p)}) \otimes_{\Z_{(p)}}
\BP_*\BP.\end{align*}

\end{lem}
\begin{proof}
For Gamma homology there is a flat-base-change result 
(see \cite[6.8 (1)]{RoWh}) which ensures that 
$$ \HG_{s,*}(\BP_*[t_i]|\BP_*; \BP_*\BP) \cong \BP_* \otimes_{\Z_{(p)}}
\HG_{s,*}(\Z_{(p)}[t_i]|\Z_{(p)}; \BP_*\BP). $$

Gamma homology of a polynomial generator is ignorant of the taken
coefficients: the Steenrod splitting of \cite[4.1]{RRo} identifies 
$\HG_{s,*}(\Z_{(p)}[t_i]|\Z_{(p)}; \BP_*\BP)$ with Gamma homology with
trivial coefficients induced up to $\BP_*\BP$ 
$$ \HG_{s,*}(\Z_{(p)}[t_i]|\Z_{(p)}; \BP_*\BP)\cong 
\HG_{s,*}(\Z_{(p)}[t_i]|\Z_{(p)}; \Z_{(p)}) \otimes \BP_*\BP.$$

In \cite{RRo} we proved this in the ungraded case. Gamma homology of a
graded algebra has homogeneous components and the above isomorphism has
no reason to preserve the internal grading. 
\end{proof}

Gamma homology of a polynomial algebra in a single variable with
coefficients in the ground ring was identified in \cite[Proposition
3.2]{RRo}. 
\begin{lem}
Summing over all internal degrees, we can identify Gamma homology of
$\Z_{(p)}[t_i]$ as 
$$ \bigoplus_{t \geq 0} \HG_{s,t}(\Z_{(p)}[t_i]|\Z_{(p)}; \Z_{(p)}) \cong
{(H\Z_{(p)})}_sH\Z.$$
\end{lem}
\begin{proof}
Using \cite{RRo} and the identification of Gamma homology with stable
homotopy groups of $\Gamma$--modules from Pirashvili--Richter \cite{PR}
we obtain
$$\HG_{s,t}(\Z_{(p)}[t_i]|\Z_{(p)}; \Z_{(p)}) \cong \left\{ 
\begin{array}{ll} 
\pi_s^{st}(\mathcal{L}_j(\Z_{(p)}[t_i]|\Z_{(p)}; \Z_{(p)})) & \mathrm{
  if } \quad j(2p^i-2) = t\\
 0 & \mathrm{otherwise}.
\end{array} \right.
$$
where $\mathcal{L}_j(\Z_{(p)}[t_i]|\Z_{(p)}; \Z_{(p)})$ is the $j^{\rm th}$
homogeneous component of  an
appropriate $\Gamma$--module. Taking all degrees together gives the
claim because 
% $\pi_s^{st}(\mathcal{L}(\Z_{(p)}[t_i]|\Z_{(p)};
% \Z_{(p)}))$ is the $s^{\rm th}$ stable homotopy of the infinite symmetric product 
% $\mathrm{SP}^\infty(\Z_{(p)})$ 
% with coefficients in the ground ring $\Z_{(p)}$ and the Dold--Thom
% theorem in turn gives us that $\pi^{st}_s(\mathrm{SP}^\infty(\Z_{(p)}))$  is
% isomorphic to 
in \cite[Theorem 4.1]{RRo} the  homotopy group 
$\pi_s^{st}(\mathcal{L}(\Z_{(p)}[t_i]|\Z_{(p)}; 
\Z_{(p)}))$ has been identified with ${(H\Z_{(p)})}_sH\Z$.
\end{proof}

In positive degrees ${(H\Z_{(p)})}_*H\Z$ is just the $p$--torsion
part of $H\Z_*H\Z$. Kochman computed this explicitly in
\cite[page 44]{K}. 

\subsection{Kochman's description of ${(H\Z_{(p)})}_*H\Z$}
Kochman provides in \cite[Theorem 3.5 (c)]{K} an explicit basis of
the $p$--torsion in ${H\Z}_*H\Z$. The result is: 
\begin{itemize}
\item 
There is only simple $p$--torsion.
\item
An explicit basis of ${(H\Z_{(p)})}_*H\Z$ over $\Z/p\Z$ consists of
all expressions 

$$ P(n_1,\ldots,n_t)\bar{\zeta}_1^{e_1}\cdot \ldots \cdot \bar{\zeta}_s^{e_s}$$
where $t \geq 0$, $t \neq 1$, $0 < n_1 < \ldots < n_t$, $e_i \geq 0$, $t+e_1+
\ldots + e_s > 0$ and $e_i = 0$ for $i < n_1$.   
Here, the degree of the $ P(n_1,\ldots,n_t)$ is $2(p^{n_1}+\ldots
+p^{n_t}) - t - 1$ and the degree of $\bar{\zeta_i}$ is $2(p^i -1)$
with the convention that the degree of $P() = 1$ is zero for $t=0$. 
\end{itemize}
Important for us are the cases $t=0$ and $t=2$. For $t=0$ the
condition $e_i = 0$ for $i < n_1$ is void, therefore elements like 
$\bar{\zeta}_1^{e_1}\ldots \bar{\zeta}_s^{e_s}$ arise with at least
one $e_i$ being positive. These elements have total degree 
\begin{equation} \label{degreest=0}
\mathrm{degree}(\bar{\zeta}_1^{e_1}\cdot \ldots\cdot \bar{\zeta}_s^{e_s}) =
\sum_{i=1}^s e_i(2p^i - 2).  
\end{equation}
For $t=2$ the situation is a little bit more involved because
non-trivial factors like $P(n,m)$ occur. The element of lowest
possible degree in this case is $P(1,2)$ with 
\begin{equation} \label{degreet=2}
  \mathrm{degree}(P(1,2)) = 2p + 2p^2 - 2 - 1 = 2p^2 + 2p - 3.
\end{equation}

\subsection[Proof of Theorem \ref{thm:bp}]{Proof of \fullref{thm:bp}}
\label{subsec:theproof}
Obstructions for extending an $n$--stage to an $(n+1)$--stage structure
live in $\HG^{n,2-n}$. From our previous arguments we know that
${(H\Z_{(p)})}_*H\Z$ corresponds to\break $\HG_{*,*}(\Z_{(p)}[t_i]|\Z_{(p)};
\Z_{(p)})$ and consists
only of simple $p$--torsion, \ie, summands of $\Z/p\Z$ as we are
working in $\Z_{(p)}$--modules. For Gamma homology of $BP_*BP$, which
is given by 
$$\bigoplus_{i\geq 0} \BP_* \otimes_{\Z_{(p)}} 
\HG_{*,*}(\Z_{(p)}[t_i]|\Z_{(p)}; \Z_{(p)}) \otimes_{\Z_{(p)}}
\BP_*\BP$$ 
we can therefore find a free $BP_*BP$--resolution of length one given
by direct sums of shifted copies of 
$BP_*BP \stackrel{p}{\lra} BP_*BP$. 

We therefore obtain that  it suffices to consider
$\mathrm{Ext}^{0,*}$-- and $\mathrm{Ext}^{1,*}$--terms in the Universal
Coefficient spectral sequence \ref{UCspsq}. 
% Therefore Gamma homology of homological degree $n$ and $n-1$ might
% occur as possible contribution to obstruction groups.
We know as well,
that the internal degree can only be of the form 
$\sum_{i=1}^N \lambda_i(2p^i -2)$; consequently possible values for
$n$ have to be of the form $\sum_{i=1}^N \lambda_i(2p^i -2) +2$ with
the $\lambda_i$ being non-negative integers. 

We are only concerned with positive homological degrees, therefore the
corresponding Gamma homology groups will be torsion. As there are no
non-trivial homomorphism from  $\BP_*\BP/p$ to $\BP_*$ the groups 
$$ \mathrm{Ext}^{0,n}_{\BP_*\BP}(\bigoplus_{i\geq 0} \BP_* \otimes_{\Z_{(p)}} 
\HG_{*,*}(\Z_{(p)}[t_i]|\Z_{(p)}; \Z_{(p)}) \otimes_{\Z_{(p)}} \BP_*\BP , \BP_*)
$$
all vanish. 

Thus we are looking for the smallest possible degree of a
non-vanishing 
$$ \mathrm{Ext}^{1,n-1}_{\BP_*\BP}(\bigoplus_{i \geq 0} \BP_* \otimes_{\Z_{(p)}} 
\HG_{*,*}(\Z_{(p)}[t_i]|\Z_{(p)}; \Z_{(p)}) \otimes_{\Z_{(p)}} \BP_*\BP
, \BP_*). $$
As $\BP_*$ is concentrated in non-negative degrees, a non-trivial map which
lowers degree by $n-1$ can come from Gamma homology in homological
degree $n -1 + \sum_{k=1}^{N'} \mu_k(2p^k -2)$. Using the constraint
for $n$ gotten above this degree can hit non-trivial groups only for 
$$ n -1 + \sum_{k=1}^{N'} \mu_k(2p^k -2) = \sum_{i=1}^N \lambda_i(2p^i
-2) + 1 + \sum_{k=1}^{N'} \mu_k(2p^k -2). $$
Considering the degrees of Kochman's generators 
$P(n_1,\ldots,n_t)\bar{\zeta}_1^{e_1}\cdot \ldots \cdot
\bar{\zeta}_s^{e_s}$ in the case $t=0$ gives the equation
$$ \sum_{i=1}^N \lambda_i(2p^i -2) + 1 + \sum_{k=1}^{N'} \mu_k(2p^k -2)
=  \sum_{j=1}^M e_j(2p^j -2),$$
so $1$ had to be even and we all believe that this is a contradiction. 

The next lowest degrees are to be found for $t=2$. We content
ourselves with considering a possible cohomology class of degree
$n-1$. Here the
corresponding equation of degrees that has to be satisfied is
$$ \sum_{i=1}^N \lambda_i(2p^i -2) + 1 = 2p^n + 2p^m -3 
+ \sum_{j=1}^M e_j(2p^j -2).$$
The generator $P(1,2)$ is of lowest possible degree and turns this
requirement into 
$$ n-1 = \sum_{i=1}^N \lambda_i(2p^i -2) + 1 = 2p + 2p^2 -3 = 2p -2 + 2p^2
-2 + 1.$$
Therefore such a homology class could occur for $n = 2p^2 + 2p -2$. \qed
\begin{rem}
Of course, the bound in \fullref{thm:bp} is just an estimate. As
some of the identification did not preserve the internal degree, this 
obstruction group may be zero. Even if it is non-trivial, it does not
have to contain an actual obstruction. Many people
believe that $\BP$ is in fact an $E_\infty$ ring spectrum, and our proof
does not disprove that. 
\end{rem}

The vanishing of Gamma cohomology for degree reasons in a certain range has a
second important consequence. Even if there is a homotopy
commutative multiplication that obeys higher coherences, it is an
issue whether this structure is unique
or whether there are different ways to find homotopies which care for
the higher coherences. Note, that the $\mathit{MU}_*$--algebra structure of
$\BP_*$ gives $\BP$ a unique $\mathit{MU}$--ring
structure \cite[2.21]{HS} and this in turn gives rise to a homotopy commutative
and associative ring structure on $\BP$, ie, a $3$--stage structure for
$\BP$. We get the following estimate.
\begin{cor}
Any $3$--stage structure for $\BP$ extends uniquely to a $(2p-1)$--stage
structure. 
\end{cor}
\begin{proof}
Obstruction groups for the unique extension of an $n$--stage structure
to an $(n+1)$--stage structure  occur as 
$ \HG^{n,1-n}(\BP_*\BP | \BP_*; \BP_*)$. The internal degree $1-n$ has to
be of the form $-\sum_{i=1}^M \lambda_i(2p^i -2)$, therefore we obtain
the constraint $n = \sum_{i=1}^M \lambda_i(2p^i -2) +1$. For the same
reasons as before homological degree $n$ corresponds to an
$\mathrm{Ext}^{0,*}$--term and vanishes. So we have to obtain an
estimate for the lowest possible $(n-1)^{\rm st}$ Gamma homology group. 
Here, the cases for $t=0$ turn out to be relevant. The equation for a
possible degree is 
$$n-1 = \sum_{i=1}^M \lambda_i(2p^i -2) = \sum_{i=1}^N e_i(2p^i -2). $$
This can occur for $e_1 = 1$ and $e_i=0$ for $i>0$, so we obtain $n-1
= 2p-2$. Consequently $\HG^{2p-1,2-2p}(\BP_*\BP | \BP_*; \BP_*)$ might be
non-trivial, and hence a given $(2p-1)$--stage structure with a fixed  
$(2p-2)$--stage structure might be prolonged to a $2p$--stage structure,
but in probably different ways. 
\end{proof}

\begin{rem}
Note, that this estimate for uniqueness is not better than the one we
would have achieved by a mere degree counting argument. We will obtain
similar results later for the localized Johnson--Wilson spectra 
$\mathcal{E}(n)$. 
\end{rem}
\section{Dyer--Lashof operations} \label{sec:dl}
An $n$--stage structure on a spectrum $E$ gives rise to some low-degree
Dyer--Lashof operations. 
If $E$ were an $H_\infty$--spectrum (see \cite{Hinfty})
then the $i^{\rm th}$ Dyer--Lashof
operation  $Q_i$ is defined as follows: take the
standard  resolution, 
traditionally called $W$ in that context, of the
cyclic group of order $p$, $\Z/p\Z$, and  assume without loss of generality
that your spectrum comes
with a $CW$--structure which is compatible with the $H_\infty$--structure. Let
$C_*$ denote cellular chain functor for spectra. Then we obtain a map
$\vartheta$ as follows. 
$$ \vartheta\co W \otimes_{\Z/p\Z} C_*(E)^{\otimes p} \ra C_*(E\Sigma_p)
\otimes_{\Sigma_p} C_*(E)^{\otimes p} \cong C_*(E\Sigma_p
\ltimes_{\Sigma_p} E^{\wedge p}) \ra C_*(E). $$
In this situation one can define $Q_i(x)$ for $x \in {H\F_p}_*(E)$ as 
$$ Q_i(x) := \vartheta_*(e_i \otimes x^{\otimes p}). $$
Here $e_i$ is the generator of $W_i$. Note that this just uses the
$i$--skeleton of $E\Sigma_p$; but this is part of an $n$--stage
structure on $E$ for $i \leq n-p$. Therefore we get the following result.
\begin{prop}
If $E$ has an $n$--stage structure with $n > p$ then there are
Dyer--Lashof operations $Q_i$ on the $\F_p$--homology of $E$ for $i \leq n-p$. 
\end{prop}
The usual Dyer--Lashof operations $Q^i$ are then given by regrading as 
$$ Q^i(x) = 
\left\{
 \begin{array}{lll}
0 & \mathrm{if} \quad \phantom{2}i < |x| \quad &\mathrm{ and }\quad  p=2 \\
Q_{i-|x|}(x) & \mathrm{if} \quad \phantom{2}i \geq |x| \quad &\mathrm{ and }\quad  p=2 \\
0 & \mathrm{if} \quad 2i < |x| \quad &\mathrm{ and } \quad p > 2 \\
\pm Q_{(2i-|x|)(p-1)}(x)\quad & \mathrm{if} \quad 2i \geq |x| \quad
&\mathrm{ and } \quad  p>2. 
 \end{array}
\right.$$

Andy Baker noticed a curious fact about Dyer--Lashof operations on
$\BP$. First, let $p$ be an odd prime.   The indecomposable element $a_{p-1} 
\in {(H\F_p)}_{2p-2}(\mathit{MU})$ is known to be in the image of  
${(H\F_p)}_{2p-2}(\BP)$. Consider an element $x = x_{2p-2}$ in\break
${(H\F_p)}_{2p-2}(\BP)$ with image $a_{p-1}$. For
such an $x$ the highest Dyer--Lashof operation $Q^i$ which we get out of the
$(2p^2 + 2p -2)$--stage structure is $Q^{2p}$. In \cite[2.11]{HKM} Hu, Kriz
and May proved that the inclusion from $\BP$ to $\mathit{MU}$ cannot
be a map  of 
commutative $S$--algebras, and they used this particular Dyer--Lashof
operation to show that (compare the correction of the proof of
\cite[2.11]{HKM} in \cite[Appendix B]{BM}). The image of $a_{p-1}$ under $Q^{2p}$ 
is $a_{(2p+1)(p-1)}$ up to decomposable elements, but there is no indecomposable 
element in ${(H\F_p)}_{(2p+1)(p-1)}(\BP)$. For $p=2$ a similar argument works 
using $a_1$. 

In the following $\SF$ denotes the colimit $\SF = \colim \SF_n$ 
where $\SF_n$ is the monoid of based homotopy equivalences of
$\mathbb{S}^n$ of degree one, so $\BSF$  is the classifying space of
spherical fibrations. 
We emphasize that the following result
is not optimal. In his thesis Lewis \cite[pages 145--6]{Lewis} wrote down
a  sketch of an argument due to
Priddy, that $\BP$ cannot be the $p$--localization of a Thom spectrum
associated to a map of $H$--spaces from a double loop space $X$ to
$\BSF$. Priddy's argument involves the Eilenberg--Moore spectral
sequence and calculations with secondary cohomology operations. 

I thank Stewart Priddy and Yuli Rudyak for
telling me that Lewis' thesis still is the only written account of
that argument. We offer our weaker result here because we think that
the proof via Dyer--Lashof operations which we give here is short and 
straightforward.
\begin{thm} \label{thm:nothom}
The Brown--Peterson spectrum $\BP$ cannot be the $p$--localization of a 
Thom spectrum associated to 
a $4$--fold loop map to $\BSF$ at $p=2$ resp.\ a $(2p+4)$--fold loop map
to $\BSF$ at any odd prime $p$.  
\end{thm}
\begin{proof}
Assume there were such a map from an $n$--fold loop space $X$ to $\BSF$ 
$$ \gamma\co X \lra \BSF$$
which would allow to write $\BP$ as the Thom spectrum associated to $\gamma$, 
$\BP = X^\gamma$. Lewis' result states that 
Thom spectra associated to $n$--fold loop maps to
$\BSF$ are $E_n$--spectra \cite[Theorem IX.7.1]{LMS}. Here $E_n$ is the
product of the little-$n$--cubes operad and the linear isometries
operad. 

The Thom isomorphism tells us that the homology of $\BP$ 
is isomorphic to the homology of $X$, and the latter maps to the homology 
of $\BSF$. The isomorphism respects the Dyer--Lashof operations
\cite[Proposition IX.7.4 (i)]{LMS}. 
In the following $H$ will denote homology with $\F_p$--coefficients. 

If $p=2$, the homology of $\BSF$ is 
\begin{equation} \label{eq:bsf}
H_*(\BSF) \cong H_*(\BSO) \otimes C_*,\end{equation}
whereas at odd primes, the homology of $\BSF$ is isomorphic to 
$$ H_*(\BSF) \cong H_*W \otimes C'_*.$$
Explicit formul{\ae} for $C_*$ and $C'_*$ can be found in \cite[page 114]{CLM}. 
The map from $\BSO$ to $\BSF$ is an infinite loop map and it is this map 
which includes the tensor factor $H_*(BSO)$ into $H_*(BSF)$. Therefore the 
tensor factor $H_*(\BSO)$ is closed under the Dyer--Lashof operations. 
A similar remark 
applies to $W$ which is a summand of $BO$ at odd primes, because there
is a splitting of infinite loop spaces  
$$ BO_{(p)} \simeq W \times W^{\perp}$$
for any odd prime $p$.  
In particular, for $x = x_{2p-2}$ in $H_*(\BP)$ with $\mathcal{P}^1_*(x) =1$ 
we obtain a non-trivial class of degree $2p-2$ in $H_*(\BSF)$. 
From \cite[pages 114--5]{CLM} it is clear that there is no 
class of that degree in the $C_*$-- resp.\ $C'_*$--part of the above tensor 
product, and therefore $x$ has to have an image in $H_*(\BSO)$ resp.\ $H_*(W)$. 
In both cases, $x$ has to hit an indecomposable element, whose $Q^{2p}$--image 
gives a generator up to decomposable elements. The lack of indecomposables in 
$H_*(\BP)$ yields a contradiction. 
\end{proof}

\begin{example}
One can ask whether our approach yields any new Dyer--Lashof operations for the
topological Hochschild homology of the Brown--Peterson spectrum. An
$E_3$--structure on $\thhbp$ justifies the calculation of $\thhbp$ given
in \cite[pages 23--4]{MS} as 
$$ \thhbp_* \cong \BP_* \otimes \Lambda(\lambda_1,\lambda_2,\ldots),
\quad \mathrm{degree}(\lambda_i) = 2p^i-1.$$
See \cite[Section 6]{BrR} as well. 

A pure degree estimate or even a more detailed analysis than in the
proof of \fullref{thm:bp} gives a  possible
obstruction class in cohomological degree $\HG^{2p,2-2p}$ caused by 
$\lambda_1$. Therefore $\thhbp$ has at least a $2p$--stage
structure for all primes. In terms of Dyer--Lashof operations this
gives  $Q^i(x)$ for
$(2i-|x|)(p-1) \leq p$, thus $2i-|x| < 2$ at odd primes and for  
$(i-|x|) \leq 2$ at $2$. 
An $E_3$--structure would provide operations $Q^i$ on $x$ for 
$2i-|x| < 2$ at odd primes and $i-|x| < 2$ at $2$ 
(compare \cite[III.3.1]{Hinfty}). 
\end{example}
\subsection*{Browder operations}
The geometry of $n$--stage structures gives rise to additional
homology operations, which give homological obstructions to extending
an $n$--stage structure to an $(n+1)$--stage structure. 
Recall that an $n$--stage structure includes actions by operad
subspaces which come from operadic composition of lower stages. For
instance an operad product 
$$ \gamma \co (E\Sigma_{2}^{(0)} \times T_{2}) \times (E\Sigma_{i}^{(0)}
\times T_{i}) \times (E\Sigma_{j}^{(0)} \times T_{j}) \rightarrow 
(E\Sigma_{i+j-1}^{(0)} \times T_{i+j-1})$$
gives elements in $\Sigma_{n+1} \times \partial T_{n+1}$ for
$i+j-1 = n+1$  which act on $E$ if the ring spectrum $E$ admits an $n$--stage
structure where $\partial T_{n+1} \subset T_{n+1}$ is the subspace of 
fully-grown $(n+1)$--trees. We already
used the identification of the homology of the quotient $T_{n+1}/\partial
T_{n+1}$ as the dual of the Lie representation. As $T_{n+1}$ is
contractible, the $H\F_p$--homology of $\partial T_{n+1}$ is isomorphic to 
${(H\F_p)}_{*+1}(T_{n+1}/\partial T_{n+1})$. On the homological level
we obtain therefore an action  of this shifted  copy of $\lie(n+1)^*$:
\begin{multline} \label{eqn:browder}
\Theta_*\co {(H\F_p)}_*(\Sigma_{n+1} \times \partial T_{n+1})
\otimes_{\Sigma_{n+1}} {(H\F_p)}_*(E)^{\otimes n+1} \\
\cong {(H\F_p)}_*(\partial T_{n+1}) \otimes {(H\F_p)}_*(E)^{\otimes n+1}
 \lra {(H\F_p)}_*(E). 
\end{multline}
If the $n$--stage structure on $E$ could be refined to an $(n+1)$--stage
structure, then this operation has to be trivial, because it factors
over  the entire tree space $T_{n+1}$ which is contractible. 

Cohen's  Theorem 12.3 in \cite{CLM} expresses
the action of the top-dimensional class of the $k$--configurations in
$\mathbb{R}^{m}$, $F(\mathbb{R}^{m},k)$, with an iterated Browder 
operation. The homology of
configuration spaces in turn is related to free Lie algebras. In
\cite[pages 263--4]{CLM} resp.\ \cite[Theorem 6.1]{Cohen:easy} Cohen
describes the homology of $F(\mathbb{R}^{m},k)$. For an even $m > 0$, 
$$ H_{(k-1)(m-1)}(F(\mathbb{R}^{m},k); \Z) \cong \lie(k) $$
as a $\Z[\Sigma_k]$--module \footnote{Note that Cohen considers
  $\lie(k)$ without the sign-representation \cite[page 32]{Cohen:easy}
  whereas  Robinson
  twists the $\Sigma_k$--action on $\lie(k)$ by the sign
  \cite[page 333]{Ro:obstr}. Therefore we
  have to adjust Cohen's statement to the above form.}. 
Up to suspension and a process of  dualization, this is precisely the
homology of the space of fully grown $(n+1)$--trees, for $k = n+1$. 

In the context of iterated loop
spaces, a non-trivial Browder operation $\lambda_n$ (coming from an
$(n+1)$--fold loop space structure on a space) is an obstruction to
extending this structure to an
$(n+2)$--fold loop space structure. The geometric operation from
\ref{eqn:browder} relates this to obstructions for $(n+1)$--stage
structures. 

It would be of interest to clarify some of the open questions about
possible $E_\infty$--structures on spectra by finding non-trivial
operations as in \ref{eqn:browder}. We hope that someone who is more
skillful with calculations than the author might actually succeed to
find obstruction classes. 

\section{Other examples} \label{sec:examples}
In the following we discuss examples where we obtain estimates for
homotopy coherence by a mere degree count, either because the 
actual (co)homology calculation would be hard if not impossible or
because we come across  a cohomology
class for which we have good reasons to believe that it is an actual
obstruction class. The localized Johnson--Wilson spectra are examples of the
first type whereas the Morava-$K$--theories are examples of the second
kind. 
\subsection{Localized Johnson--Wilson $\mathcal{E}(i)$}
We start by considering the localized version of the John\-son--Wilson
spectra $E(i)$, $\mathcal{E}(i)$ with  $\mathcal{E}(i)_* 
\cong {(\bpi_*)}_{I_i}$ (compare \cite{B}). 
The coefficients of $\mathcal{E}(i)_*$ are the coefficients of $\bpi$
localized away from the ideal $I_i = (p,v_1,\ldots,$
$v_{i-1})$. Therefore units in
$\mathcal{E}(i)_*$ are of the form $s v_i^j + r$ with $s$ a unit in
the $p$--local integers, $j \in \Z$ and $r \in I_i$ (compare \cite[Section
1]{B}). 
By \cite{B} it is known that   
$\mathcal{E}(i)_*\mathcal{E}(i)$ is free over $\mathcal{E}(i)_*$. Thus 
possible obstructions live in 
$$\mathrm{Hom}^{2-n}_{\mathcal{E}(i)_*}(\lie(m)^* \otimes
\mathcal{E}(i)_*[\Sigma_m^{n-m+1}] \otimes
\mathcal{E}(i)_*\mathcal{E}(i)^{\otimes m}, \mathcal{E}(i)_*). $$
As $\mathcal{E}(i)$ is Landweber exact over $\BP$, the algebra of
cooperations can be described as
$$ \mathcal{E}(i)_*\mathcal{E}(i) = \mathcal{E}(i)_* \otimes_{\BP_*}
\BP_*\BP \otimes_{\BP_*} \mathcal{E}(i)_*.$$
The degrees of this algebra and hence the degrees of the
$\mathcal{E}(i)_*$--module 
$$\lie(m)^* \otimes \mathcal{E}(i)_*[\Sigma_m^{n-m+1}] \otimes
\mathcal{E}(i)_*\mathcal{E}(i)^{\otimes m} $$
are concentrated in degrees of the form 
\begin{equation}\label{jwdegrees} 
\sum_{\genfrac{}{}{0pt}{}{i=1}{i \neq n}}^M \lambda_i (2p^i -2) + 
\mu_n(2p^n -2)
\end{equation}
with the  $\lambda_i$ being non-negative integers and $\mu_n \in \Z$. 

As this is the underlying module for the chain complex whose dual
gives Gamma cohomology we obtain the following result. 
\begin{prop}\label{prop:jw}
The $i^{\rm th}$ localized Johnson--Wilson spectrum $\mathcal{E}(i)$
possesses at least a $2p$--stage structure which is unique up to the 
$(2p-1)$--stage.
\end{prop}
\begin{proof}
For any potential  obstruction we need a non-trivial
$\mathcal{E}(i)_*$--linear morphism  which raises degree by some number
of the form $n-2$ with $n \geq 3$. As the domain  
$\lie(m)^* \otimes \mathcal{E}(i)_*[\Sigma_m^{n-m+1}] \otimes
\mathcal{E}(i)_* \mathcal{E}(i)^{\otimes m} $ and the target
$\mathcal{E}(i)_*$ are concentrated in degrees as calculated in 
\ref{jwdegrees} the minimal such difference of degrees is
$2p$. Therefore there might be an obstruction to extending a given
$2p$--stage structure to a $(2p+1)$--stage structure. 
\end{proof}

\subsection{Ordinary Johnson--Wilson spectra}
The $i^{\rm th}$ Johnson--Wilson spectrum $E(i)$ is Land\-weber exact and a
homotopy commutative ring spectrum for
all primes. The algebra of cooperations $E(i)_*E(i)$ is flat over
$E(i)_*$ and  hence the  K\"unneth
theorem applies, but we have to use the universal coefficient
spectral sequence to control the source 
$E(i)^1(Q_{n+1}^m \ltimes_{\Sigma_m} E(i)^{\wedge m})$ of possible
obstructions for highly coherent homotopy commutative structures. The
$E_2$--term of this spectral sequence is as follows:
\begin{multline*}
\mathrm{Ext}_{E(i)_*}^{p,q}(E(i)_*(Q_{n+1}^m \ltimes_{\Sigma_m} 
E(i)^{\wedge m}), E(i)_*) \\
\cong 
\mathrm{Ext}_{E(i)_*}^{p,q}(\Sigma^{n-1} \lie(m)^* \otimes
E(i)_*[\Sigma_m]^{\otimes (n-m+1)} \otimes E(i)_*E(i)^{\otimes m}, E(i)_*)
\end{multline*}
Here the $(n-1)^{\rm st}$ suspension in $E(i)$--homology comes from the
geometry of the filtration quotient $Q_{n+1}^m$ as it was discussed in 
\fullref{sec:geom}. 

In this spectral sequence we just have to consider the $0$-- and 
$1$--line, because\break $E(i)_*E(i)$ is a countable colimit of free 
$E(i)_*$--modules and the other tensor factors are free. As we are only
interested in cohomological degree one and
$\mathrm{Ext}_{E(i)_*}^{p,q}$ corresponds to cohomological degree
$p+q$ we have to determine, for which $n \geq 3$ there are non-trivial 
$\mathrm{Ext}_{E(i)_*}^{0,1}$-- and $\mathrm{Ext}_{E(i)_*}^{1,0}$--
terms. Taking the internal $(n-1)$--shift into account we have to look
out for homomorphisms which raise degree by $n-2$ resp.\ $n-1$. Counting
degrees once again we obtain the following result.
\begin{prop}
The $i^{\rm th}$ Johnson--Wilson spectrum $E(i)$ possesses at least a
$(2p-1)$--stage structure which is unique up to the $(2p-2)$--stage.
\end{prop}
\begin{proof}
Both sides are concentrated in degrees of the form 
$$\sum_{\genfrac{}{}{0pt}{}{i=1}{i \neq n}}^M \lambda_i (2p^i -2) + 
\mu_n(2p^n -2).$$
Due to the possibility of a non-trivial $\mathrm{Ext}^{1,1-n}$--term,
the lowest possible value for $n$ is 
$n = 2p - 2 +1 = 2p - 1$.
\end{proof}

\vspace{-1mm}
\begin{rem}
Note that the above estimate is much too weak for $i=1$ because $E(1)$
possesses a unique $E_\infty$--structure at all primes \cite[Theorem 6.2]{BR}.
\end{rem}

\vspace{-1mm}
\subsection{Morava-$K$--theory}

\vspace{-1mm}
Last but not least we will close with a negative result. 
It is known that none of the Morava-$K$--theories could possess 
$E_\infty$--structures: if $K(n)$ had one, then its connective cover $k(n)$
had one as well; in particular, $k(n)$ would be an $H_\infty$--spectrum
which then had to split as a wedge of suspensions of Eilenberg--MacLane
spectra (see \cite[III, Theorem 4.1]{Hinfty}).  The argument in 
\cite[III, Theorem 4.1]{Hinfty} even shows that $k(n)$
cannot possess any $H_2$--structure. From \cite[2.10]{M} it follows that
the  connected cover functor $c$ \cite[VII.3.2]{M} sends $H_2$--spectra
to  $H_2$--spectra. Therefore $K(n)$ cannot have an $E_2$--structure,
because  this would give
rise to such an $H_2$--structure on $k(n)$.

At $p=2$, $K(n)$ is even not homotopy
commutative. For $\mathit{MU\!}_{(p)}$--algebra structures on $K(n)$ see
\cite{G,L}. At odd primes, the $K(n)$ are homotopy
commutative and associative ring spectra. 
 
Consider $K(n)$ for an odd prime $p$. We know that the algebra of
cooperations\break $K(n)_*(K(n))$ for $K(n)$ consists of an \'etale part
tensored with an exterior algebra 
$\Lambda_{K(n)_*}(\tau_0,\tau_1,\ldots,\tau_{n-1})$ with the $\tau_i$
being elements of degree $2p^i-1$. This non-\'etale  part could  give rise to
possibly non-trivial obstructions. Here $\tau_0$ is of internal degree 
one. A possible non-trivial obstruction groups could therefore be
obtained  by a class in cohomological degree $n = 2+1 =3$. Hence even
at odd primes there might be a non-trivial
obstruction to extending the homotopy commutative structure on $K(n)$
to a $4$--stage structure.

\vspace{-1mm}
\begin{rem}
The situation for the spectrum $P(n)$ at odd primes with $P(n)_* =
\F_p[v_n,v_{n+1}, \ldots]$ is the same as the
one for Morava-$K$--theory. This spectrum cannot possess an
$H_2$--structure. The algebra of cooperations is 
$$ P(n)_*P(n) = P(n)_*\otimes_{\BP_*} \BP_*\BP \otimes \Lambda(a_0,
\ldots,a_{n-1}) $$
with $a_i$ being of degree $2p^i-1$; in particular $a_0$ has degree
one, which should give rise to an obstruction to a $4$--stage
structure. 
\end{rem}

In both cases, $K(n)$ and $P(n)$, there are candidates for non-trivial
Browder operations for $K(n)$--homology respectively $H\F_p$--homology,  
but we have no proof so far that these homology classes actually arise in this
way.

\bibliographystyle{gtart} \bibliography{link}

\end{document}